\documentclass[preprint,onefignum,onetabnum]{siamonline250211}

\usepackage{lipsum}
\usepackage{amsfonts}
\usepackage{graphicx}
\usepackage{epstopdf}
\usepackage{algorithmic}
\usepackage{siunitx}
\usepackage{microtype}
\microtypesetup{spacing=true}
\ifpdf
  \DeclareGraphicsExtensions{.eps,.pdf,.png,.jpg}
\else
  \DeclareGraphicsExtensions{.eps}
\fi

\headers{Solving MFC problems using VCNF}{J. Zhao, M. Zhou, X. Zuo, and W. Li}

\title{Variational conditional normalizing
flows for computing second-order mean field control problems}

\author{%
Jiaxi Zhao\thanks{Department of Mathematics, National University of Singapore, Singapore.  (\email{jiaxi.zhao@u.nus.edu})}%
\and 
Mo Zhou\thanks{Department of Mathematics, University of California Los Angeles, CA, 90095. 
  (\email{mozhou366@math.ucla.edu}).}%
\and
Xinzhe Zuo \thanks{Department of Mathematics, University of California Los Angeles, CA, 90095.  (\email{zxz@math.ucla.edu}).}%
\and 
Wuchen Li \thanks{Department of Mathematics, University of South Carolina, SC, 29036. (\email{wuchen@mailbox.sc.edu}).}%
}

\usepackage{amsopn}

\makeatletter
\newcommand*{\addFileDependency}[1]{% argument=file name and extension
  \typeout{(#1)}% latexmk will find this if $recorder=0 (however, in that case, it will ignore #1 if it is a .aux or .pdf file etc and it exists! if it doesn't exist, it will appear in the list of dependents regardless)
  \@addtofilelist{#1}% if you want it to appear in \listfiles, not really necessary and latexmk doesn't use this
  \IfFileExists{#1}{}{\typeout{No file #1.}}% latexmk will find this message if #1 doesn't exist (yet)
}
\makeatother

\newcommand*{\myexternaldocument}[1]{%
    \externaldocument{#1}%
    \addFileDependency{#1.tex}%
    \addFileDependency{#1.aux}%
}
\usepackage{essay-def}
\usepackage{comment}
\usepackage{amssymb}
\usepackage{subcaption}
\usepackage{cleveref}
\usepackage{booktabs}

\ifpdf
\hypersetup{
  pdftitle={Solving MFC problems using VCNF},
  pdfauthor={J. Zhao, M. Zhou, X. Zuo, and W. Li}
}
\fi

\myexternaldocument{ex_supplement}

\begin{document}

\maketitle

% REQUIRED
\begin{abstract}
Mean field control (MFC) problems have vast applications in artificial intelligence, engineering, and economics, while solving MFC problems accurately and efficiently in high-dimensional spaces remains challenging. This work introduces variational conditional normalizing flow (VCNF), a neural network-based variational algorithm for solving general MFC problems based on flow maps. Formulating MFC problems as optimal control of Fokker--Planck (FP) equations with suitable constraints and cost functionals, we use VCNF to model the Lagrangian formulation of the MFC problems. In particular, VCNF builds upon conditional normalizing flows and neural spline flows, allowing efficient calculations of the inverse push-forward maps and score functions in MFC problems. We demonstrate the effectiveness of VCNF through extensive numerical examples, including optimal transport, regularized Wasserstein proximal operators, and flow matching problems for FP equations.

\end{abstract}


% REQUIRED

\section{Introduction}

In recent years, generative artificial intelligence (AI) has attracted much attention in the scientific computing community across various fields. Typical applications include text generation \cite{achiam2023gpt,ray2023chatgpt}, image generation \cite{rombach2022high, song2020score, ho2022video}, and protein folding \cite{jumper2021highly,abramson2024accurate}. At the core of these applications are nonlinear generative models, such as normalizing flows  \cite{10.5555/3454287.3454962} and neural ordinary differential equations \cite{chen2018neural}. These models aim to learn a transport map from a reference distribution to a target distribution, enabling sample generation that approximates the target. The optimization problems in generative models can be interpreted as variational problems in probability space, also known as the mean field control (MFC) problems \cite{fornasier2014mean,E2018AMO}. 

MFC problems incorporate interactions among multiple agents in the optimal control problem, providing a method for modeling large-scale collective behaviors in artificial intelligence \cite{liu2018mean}, engineering \cite{lee2021controlling, lee2022mean}, and economics \cite{moll14,cardialiaguet2018}. In literature, the second-order MFC problems can be expressed as a coupled system of Hamilton--Jacobi--Bellman (HJB) and Fokker--Planck (FP) equations, which jointly characterize the optimal velocity field governing the agents' density evolution. The second-order MFC emphasizes Kolmogorov forward and backward operators from the diffusion processes in the coupled systems. In particular, optimal transport (OT), i.e., Wasserstein distances between probability distributions, with their dynamical formulation \cite{benamou2000computational}, can be viewed as a special class of first order MFC problem with fixed initial and terminal time constraints on distributions and without diffusion.

Traditional methods for solving OT and MFC problems are well-studied in low-dimensional spaces \cite{benamou2017variational,benamou2000computational,yu2024fast}. Recently, machine learning-based approaches have emerged as promising tools for tackling high-dimensional MFCs, pioneered by \cite{doi:10.1073/pnas.1922204117}, which solves high-dimensional deterministic MFCs by approximating the value function using deep neural networks. In Lagrangian coordinates, solving MFC problems often involves solving the Monge--Amp\'ere equation, which arises from the change of variable formula for probability densities. Classical numerical PDE algorithms \cite{benamou2000computational,froese2012numerical} perform well in low dimensions. In comparison, modern generative models tackle the Monge--Amp\'ere equation from a variational perspective, optimizing over a parametrized family of transformations by minimizing certain objective functions, such as the Kullback–Leibler (KL) divergence between the generated and target distribution. A typical example is the neural spline flow \cite{10.5555/3454287.3454962}, which constructs transport maps based on monotone rational-quadratic splines. This network structure enables efficient computation of both the inverse transport map and the determinant of the Jacobian of the transport map. Such a formulation is particularly useful in approximating the score function, defined as the gradient of the log-density function. Additionally, the score function approximates the Kolmogorov forward operator associated with diffusion processes. Given the potential of generative models and score functions, a natural question arises: \textit{Can we apply generative models to solve second-order MFC problems in high dimensional spaces?}

\begin{figure}
\centering
\includegraphics[width=\linewidth]{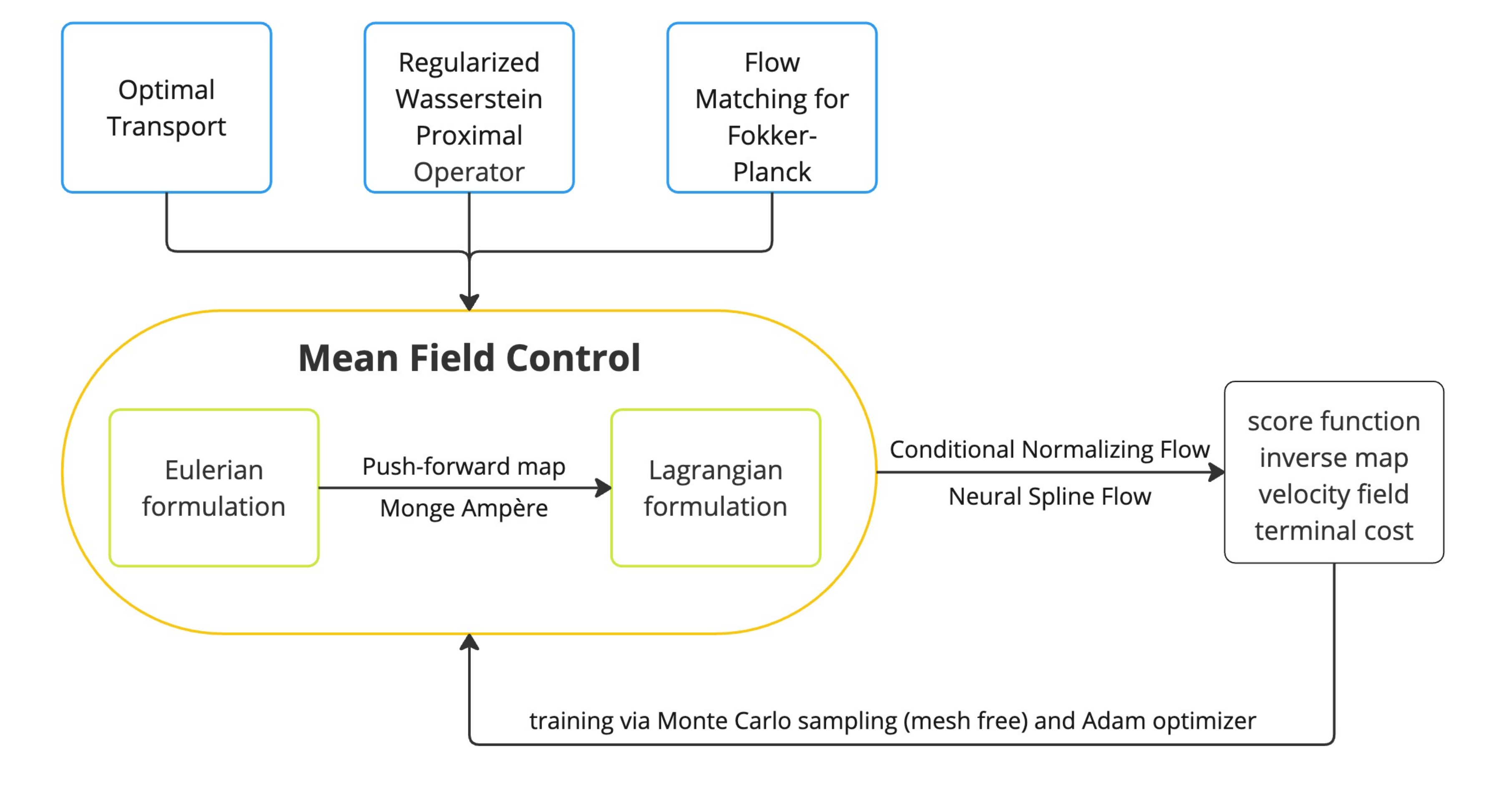}
\caption{VCNF framework for solving generalized MFC problems. }
\label{fig:flowchart}
\end{figure}

In this paper, we introduce the variational conditional normalizing flow (VCNF), an extension of conditional normalizing flows \cite{winkler2019learning} designed for solving general MFC problems. Conditioning on time $t$, VCNF models a family of distributions $p(\mfx, t)$ that evolve continuously over time to solve the MFC problems. Unlike standard normalizing flows with fixed time stamps, VCNF allows for evaluation and sampling at arbitrary time. The time derivative can be computed using finite difference or automatic differentiation. By leveraging numerical approximations of score functions, VCNF can effectively handle second-order MFC problems with diffusive density evolution. Specifically, we demonstrate this capability by solving the FP equation within the flow matching framework. We show in the numerical experiments that VCNF is highly robust and efficient in solving various MFC problems. In particular, we are able to achieve around $10^{-3}$ relative error, which is commonly expected for neural network-based models in scientific computing \cite{raissi2019physics}, with little parameter tuning. In addition, VCNF effectively captures the evolution of multimodal distributions and accurately models the transitions between unimodal and multimodal distributions under different initial conditions, demonstrating its flexibility and reliability. A schematic of our framework is demonstrated in \cref{fig:flowchart}, details of which will be discussed in \cref{sec:mfc_lagrangian,sec:mfc_examples,sec:cnf}.

This paper is organized as follows. In \cref{sec:intro to mfc}, we review some background on MFC problems, including both Eulerian and Lagrangian formulations. Several important examples are covered. In \cref{sec:model and alg}, we present the construction of the model, the details for evaluating various quantities, and the algorithm for solving MFC problems. In \cref{sec:numerics}, we provide extensive numerical experiments to show the robustness and efficiency of VCNF. In \cref{sec:discussion}, we conclude this work and discuss potential future directions. 

\section{Mean field control problems}\label{sec:intro to mfc}

In this section, we review second-order MFC problems in Eulerian and Lagrangian coordinates. We also introduce several important examples, including OT, RWPO, and flow matching for FP equation.

\subsection{Mean field control problem in Eulerian coordinates}

We consider the MFC problem where one aims to minimize the variational objective
\begin{equation}\label{eq:MFC_objective}
\inf_{\mfv,p} \int_0^T \int_{\RR^d} L(\mfx,t,\mfv(\mfx,t),p(\mfx,t)) p(\mfx,t) \,\rd \mfx\,\rd t + \int_{\RR^d} G(\mfx,p(\mfx,T)) p(\mfx,T) \,\rd \mfx\,,
\end{equation}
subject to the FP equation with initial condition
\begin{equation}\label{eq:FP}
\pt p(\mfx,t) + \nabla_\mfx\cdot(p(\mfx,t)\mfv(\mfx,t)) = \gamma \Delta_\mfx p(\mfx,t)\,, \quad\quad p(\cdot,0)=p_0\,.
\end{equation}
The objective functional involves a running cost $L:\RR^d\times[0,T]\times\RR^d\times\RR_+ \to \RR$ and a terminal cost $G:\RR^d\times\RR_+ \to \RR$. The state dynamic is characterized by the evolution of the density function \eqref{eq:FP}. In the MFC problem, a central planner controls the drift of the dynamic $\mfv(\mfx,t)$, which has a direct impact on the density $p(\mfx,t)$. The goal of the MFC problem is to minimize the cost \eqref{eq:MFC_objective} w.r.t. $\mfv$ and $p$ subject to the constraint \eqref{eq:FP}. The diffusion term $\Delta_\mfx p(\mfx,t)$ indicates that this MFC problem is of second-order.

The MFC problem also has a stochastic formulation:
$$\inf_\mfv \EE \sqbra{\int_0^T L\parentheses{\mfx_t,t,\mfv(\mfx_t,t),p(\mfx_t,t)} \,\rd t + G(\mfx_T,p(\mfx_T,T)) }\,,$$
subject to the state dynamic
$$\rd \mfx_t = \mfv(\mfx_t,t) \,\rd t + \sqrt{2\gamma} \,\rd W_t\,, \quad\quad \mfx_0 \sim p_0\,.$$
$p(\cdot,t)$ is the density function of $\mfx_t$ and satisfies the FP equation \eqref{eq:FP}. $W_t$ is a $d$-dimensional standard Brownian motion. $\mathbb{E}$ is the expectation over all realizations of the state dynamic. 

We remark that mean field games (MFGs) \cite{lasry2007mean, cardaliaguet2010notes} differ from MFCs primarily due to the absence of a central planner in MFGs. Consequently, a change of control by a single agent does not affect the population distribution. However, for a large class of MFG problems, one can find the corresponding variational formulations (also known as potential MFGs) and transform them into MFC problems \cite{doi:10.1073/pnas.1922204117,HUANG2023112155,zhang2023mean}.

Numerical methods for solving MFC problems generally fall into two broad categories. The first category  focuses on the variational formulation of MFC, optimizing the objective functional directly. In this framework, probabilistic modeling techniques, such as neural ODE~\cite{chen2018neural} and normalizing flow~\cite{kobyzev2020normalizing} are commonly used to get the population distribution. These methods offer advantages in flexibility and scalability, making them well-suited for high-dimensional problems. The second category leverages the HJB formulation of the MFC, where the optimal solution is characterized by a coupled FP-HJB system. Common approaches in this category parameterize the solution to the HJB equation and solve the FP-HJB system \cite{briceno2018proximal, liu2020computational, zhou2024deep, darbon2020overcoming, zhou2021actor}. Hybrid approaches also exist, which employ mixed loss functions to enhance robustness and accuracy \cite{doi:10.1073/pnas.2024713118,zhou2024solving}.

For completeness, we present the HJB-FP for MFC problem. Interested readers could refer to \cite[Chapter 4]{bensoussan2013mean} for more details. We define the Hamiltonian $H:\RR^d\times[0,T]\times\RR^d\times\RR_+ \to \RR$ as the Legendre transform of $L$ w.r.t. $\mfv$, 
$$H(\mfx,t,\mfw,p) := \sup_{\mfv\in\RR^d} \mfw\tp \mfv - L(\mfx,t,\mfv,p)\,.$$
Then, the optimal control, i.e., optimal velocity, is given by
$$\mfv(\mfx,t) = \argmax_{\mfv\in\RR^d} \nabla_\mfx \phi(\mfx,t)\tp \mfv - L(\mfx,t,\mfv,p(\mfx,t))\,,$$
where $p(\mfx,t)$ and $\phi(\mfx,t)$ satisfy the coupled FP-HJB system
\begin{equation*}
\left\{\begin{aligned}
& \pt p(\mfx,t) + \nabla_\mfx\cdot(p(\mfx,t) \mathrm{D}_\mfw H(\mfx,t,\nabla_\mfx \phi(\mfx,t),p(\mfx,t)) ) = \gamma \Delta_\mfx p(\mfx,t)\,, \\
& \pt \phi(\mfx,t) + H(\mfx,t,\nabla_\mfx \phi(\mfx,t),p(\mfx,t))\\
& \hspace{1.1in} - \pd{L(\mfx,t,\mfv(\mfx,t),p(\mfx,t))}{p}\cdot p(\mfx,t) = -\gamma\Delta_\mfx\phi(\mfx,t)\,,\\
& p(\cdot,0)=p_0\,, \quad \quad \phi(\mfx,T) = -G(\mfx,p(\mfx,T)) - \pd{G(\mfx,p(\mfx,T))}{p}p(\mfx,T)\,.
\end{aligned}\right.
\end{equation*}

% In~\cite{HUANG2023112155}, the authors propose to use normalizing flows (NFs) to solve variational MFGs via a trajectory regularization term to approximate the kinetic energy in the loss function. This term consists of finite differences between adjacent layers of the NFs as an approximation of the velocity. From this perspective, different layers of the NFs can be viewed as the time steps of the trajectory, which is reminiscent of the well-known relation between deep neural network and ODE~\cite{chen2018neural}.

\subsection{Mean field control problem in Lagrangian coordinates}\label{sec:mfc_lagrangian}
We can also describe the MFC problem in Lagrangian coordinates, where we focus on the push-forward map $f$ instead of the velocity field $\mfv$. Let $q(\mfz)$ be a reference measure. In this paper, we use the standard Gaussian distribution $\mathcal{N}(0,I_d)$ as the reference measure, and we will not distinguish a probability measure with its density function. Let $f:\RR^d \times [0,T] \to \RR^d$ be a flow of invertible push forward map, which gives the density function at $t\in [0,T]$ by 
$$p(\cdot, t) = f(\cdot, t)_{\#} q(\cdot)\,.$$
The associated Monge--Amp\'ere equation is
\begin{equation*}
p(\mfx, t) = p(f(\mfz,t), t) = q(\mfz) \left| \text{det} \frac{\partial f(\mfz, t)}{\partial \mfz} \right|^{-1}.
\end{equation*}
In this case, the density function $p$ satisfies the transport equation
\begin{equation*}
\pt p(\mfx,t) + \nx \cdot\sqbra{p(\mfx,t)\, \pt f( f^{-1}(\mfx,t) ,t)} = \pt p(\mfx,t) + \nx \cdot\sqbra{p(\mfx,t)\, \pt f(\mfz,t)} = 0\,,
\end{equation*}
where $f^{-1}(\cdot,t)$ is the inverse map of $f(\cdot,t)$. If we assume that $p(\mfx,t)$ satisfies the FP equation \eqref{eq:FP}, then we can apply Nelson's transformation \cite{nelson1966derivation}. Since 
\begin{equation*}
  \nabla p(\mfx, t)=p(\mfx, t)\cdot \frac{\nabla p(\mfx, t)}{p(\mfx, t)}=p(\mfx, t)\nabla\log p(\mfx, t)\,, 
\end{equation*}
where $\nabla\log p(\mfx, t)$ is often called the score function, the FP equation can be rewritten as
\begin{equation}\label{equ:fp-score}
    \pt p(\mfx,t) + \nx \cdot \sqbra{ p(\mfx,t)\parentheses{\mfv(\mfx,t) - \gamma \nx\log p(\mfx,t)} } = 0\,.
\end{equation}
We can compute the velocity field $\mfv$ through
$$\mfv(\mfx,t) = \pt f( f^{-1}(\mfx,t) ,t) + \gamma \nx\log p(\mfx,t) = \pt f(\mfz ,t) + \gamma \nx\log p(f(\mfz,t),t)\,.$$
Therefore, the MFC problem in Lagrangian coordinates is
\begin{equation}\label{eq:MFC_Lagrange}
\begin{aligned}
\inf_{f,p} \int_0^T & \int_{\RR^d} L\parentheses{f(\mfz,t),t,\pt f( \mfz ,t) + \gamma \nx\log p(f(\mfz,t),t), p(f(\mfz,t),t)} q(\mfz) \,\rd \mfz\,\rd t \\
& + \int_{\RR^d} G(f(\mfz,t),p(f(\mfz,t),T)) q(\mfz) \,\rd \mfz \\
\text{s.t.} \quad & f(\cdot,t)_\# q = p(\cdot,t) \,,\, f(\cdot,0)_\# q = p_0\,.
\end{aligned}
\end{equation}

In VCNF, the density $p(f(\mfz,t),t)$ is obtained through the logarithm of the Monge--Amp\'ere equation
\begin{equation}\label{equ:log-den}
\log p(\mfx, t) = \log p(f(\mfz,t), t) = \log q(\mfz) - \log \left| \text{det} \frac{\partial f(\mfz, t)}{\partial \mfz} \right|.
\end{equation}
% Notice that the Jacobian of the mapping $\frac{\partial f(\mfz, t)}{\partial \mfz}$ is evaluated for each layer and summed up to obtain the density for the whole model.
The temporal derivative of the push forward map $f( \mfz ,t)$ and the score function $\nx\log p(f(\mfz,t),t)$ are approximated through numerical differentiation (see \cref{sec:VCNF quantities}).

\subsection{Examples of mean field control problems}\label{sec:mfc_examples}
We provide three important examples for MFC problems in this section, including OT, RWPO, and flow matching for FP equations.

\medskip
\noindent
\textbf{Optimal transportation.} The OT problem can be viewed as a special MFC problem where the terminal cost is replaced by a terminal constraint. By the Benamou--Brenier formulation \cite{benamou2000computational}, computing the Wasserstein-$2$ distance between two probability distributions is equivalent to the following variational problem
\begin{equation}\label{eq:OT}
\begin{aligned}
&\inf_{v,p} \int_0^1 \int_{\RR^d} \half \norml\bm{v}(\mfx,t)\normr^2 p(\bm{x}, t)\, \rd\bm{x} \, \rd t\,,           \\
&\mathrm{s.t.} \quad \partial_t  p(\mfx,t) + \nabla \cdot (p(\mfx,t)\bm{v}(\mfx,t)  ) = 0\,, \quad p(\cdot, 0) =  p_0\,, \quad p(\cdot, 1) =  p_1\,.
\end{aligned}   
\end{equation}
The minimal cost to \cref{eq:OT} is $\frac12W_2(p_0,p_1)^2$, where $W_2(\cdot,\cdot)$ denotes the Wasserstein-$2$ distance between two measures. In Lagrangian coordinates, the OT problem reads
\begin{equation*}
\begin{aligned}
&\inf_{f,p} \int_0^1 \int_{\RR^d} \half \norml \pt f(\mfz,t) \normr^2 q(\mfz)\, \rd\mfz \, \rd t\,,           \\
&\mathrm{s.t.} \quad f(\cdot,0)_\# q = p_0\,, \quad f(\cdot,1)_\# q = p_1\,.
\end{aligned}   
\end{equation*}

\medskip
\noindent
\textbf{Regularized Wasserstein proximal operator.}
Closely related to OT, the RWPO computes a regularized OT problem with a diffusion term and a terminal cost \cite{li2023kernel}. The objective is
\begin{equation}\label{equ:rwpo-loss}
\begin{aligned}
&\inf_{\mfv,p}  \int_0^T \int_{\RR^d} \frac12\norm{\mfv(\mfx,t)}^2  p(\mfx,t) 
\,\rd \mfx \,\rd t + \int_{\RR^d} V(\mfx)  p(\mfx,T) \,\rd \mfx \,,\\
&\text{s.t.} \quad \pt p(\mfx,t) + \nabla_{\mfx}\cdot \parentheses{ p(\mfx,t) \mfv(\mfx,t)}
=\frac1\beta \Delta_{\mfx}  p(\mfx,t)\,, \quad p(\cdot,0)=p_0\,.
\end{aligned}
\end{equation}
For derivation of this variational formulation, see \cref{sec:rwpo}. Using Nelson's transformation introduced in \Cref{sec:mfc_lagrangian}, the RWPO problem in Lagrangian coordinate is
\begin{equation}\label{eq:RWPO}
\begin{aligned}
&\inf_{f,p}  \int_0^T \int_{\RR^d} \frac12\norm{ \pt f(\mfz,t) + \frac{1}{\beta} \nx\log p(f(\mfz,t),t)}^2  q(\mfz) 
\,\rd \mfz \,\rd t + \int_{\RR^d} V(f(\mfz,T)) \,  q(\mfz) \,\rd \mfz\,,\\
& \text{s.t.} \quad  f(\cdot,t)_\# q = p(\cdot,t) \,,\quad f(\cdot,0)_\# q = p_0\,.
\end{aligned}
\end{equation}

\medskip
\noindent
\textbf{Flow matching for Fokker--Planck equation.}
Another important application of MFC is the flow matching problem for the FP equation
$$\pt p(\mfx,t) + \nx\cdot\parentheses{p(\mfx,t)\mfb(\mfx,t)} = \gamma \Delta_\mfx p(\mfx,t)\,.$$
For example, \cite{tang2022adaptive, zeng2023adaptive}
use NFs to solve the steady state of the (fractional) FP equation. Similarly,
CNFs are used to approximate the time-dependent FP equation \cite{feng2021solving}.

In our setting, we aim to learn the push forward map $f:\RR^d\times[0,T]\to\RR^d$ such that $f(\cdot,t)_\#q=p(\cdot,t)$ represents the solution to the FP equation at $t$. The objective functional is
\begin{equation}
\begin{aligned}
&\inf_{f,p} \int_0^T \int_{\RR^d} \norm{\pt f(f(\mfz,t),t) - b(f(\mfz,t),t) + \gamma \nx\log p(f(\mfz,t),t)}^2  q(\mfz) \,\rd \mfz\,\rd t\,, \\
 &\text{s.t.} \quad  f(\cdot,t)_\# q = p(\cdot,t) \,,\quad f(\cdot,0)_\# q = p_0\,.
\end{aligned}
\end{equation}

\section{Variational conditional normalizing flows} \label{sec:model and alg}
In this section, we describe the architecture of the VCNF model based on invertible spline transformations. We then explain its applications to solving MFC problems.

\subsection{Conditional normalizing flow via neural spline flow}\label{sec:cnf}

As discussed in \cref{sec:mfc_lagrangian}, the Lagrangian formulation of MFC problems requires parameterizations of a family of invertible maps in high dimensions. Following the idea in \cite{10.5555/3454287.3454962}, we build multi-dimensional invertible transformations through one dimensional invertible transformations $\phi(\cdot, \theta): \mbR \rightarrow \mbR$ parameterized by $\theta$ based on spline functions (see details in \cref{sec:neural_spline_flow}). Comparing to models \cite{dinh2016density} with affine transformation, this class of models are powerful in fitting complex distributions \cite{10.5555/3454287.3454962}. 

The overall transformation is built on compositions of several invertible auto-regressive layers:
\begin{equation}\label{equ:auto-regressive}
    f_\psi(\cdot,t) = f^{(l)}_\psi(\cdot,t) \circ \cdots \circ f^{(1)}_\psi(\cdot,t)\,,
\end{equation}
where $\psi$ denotes the parameters of the whole model. $\{f^{(i)}_\psi(\cdot,t)\}_{i=1}^l$ are auto-regressive layers with different coupling orders. A typical auto-regressive layer with coupling order $(1,2,\ldots,d)$ consists of $d$ coupling transformations $\mfx^{(k)} \mapsto \mfx^{(k+1)}$ composed sequentially (for $k = 0,\ldots,d-1$), each of which consists of the following steps:
\begin{enumerate}
\item Compute $\theta_k = \text{NN}^{(k)}(\mfx_{1:k}^{(k)}, t) \in \RR^{3K - 1}$ where $\text{NN}^{(k)}$ is the conditioning network. 
\item Compute $x_{k+1}^{(k+1)} = \phi(x_{k+1}^{(k)}, \theta_k)$, the detailed parametrization of $\phi(x_{k+1}^{(k)}, \theta_k)$ is explained in \cref{sec:neural_spline_flow}.
\item Set $x_i^{(k+1)} = x_i^{(k)}$ for $i = 1, ..., k, k+2, ..., d$ and return $\mfx^{(k+1)}$.
\end{enumerate}
Note that $\theta_k$, the parameters of one dimensional invertible transformations, are not the \emph{inputs} but the \emph{outputs} of the conditioning networks $\text{NN}^{(k)}: \RR^{k+1} \to \RR^{3K - 1}$, which are fully-connected feedforward neural networks for $k=0,\ldots,d-1$. The $K$ is the number of the bins for the spline functions. For example, an auto-regressive layer of the CNF in $2$ dimensions with coupling order $(1,2)$ can be written as
\begin{equation}\label{equ:cnf-auto-regressive}
\begin{pmatrix}
x_1^{(0)} \\
x_2^{(0)}
\end{pmatrix} 
\underset{x_2^{(1)}=\ x_2^{(0)}}{\xrightarrow{x_1^{(1)}=\ \phi(x_1^{(0)}, \text{NN}^{(0)}_1(t))}}
\begin{pmatrix}
x_1^{(1)} \\
x_2^{(1)}
\end{pmatrix} 
\underset{x_2^{(2)}=\ \phi(x_2^{(1)},\, \text{NN}^{(1)}_2(x_1^{(1)},t))}{\xrightarrow{\quad\quad x_1^{(2)} =\ x_1^{(1)} \quad\quad }}
\begin{pmatrix}
x_1^{(2)} \\
x_2^{(2)}
\end{pmatrix}.
\end{equation}
Specifically, the 2D transformation is decomposed into two steps: in the first step, the first coordinate performs invertible transformations parametrized by the output of the conditioning network $\text{NN}^{(0)}$ and the second step contains an invertible transformation of the second coordinate. 

As illustrated in \cref{equ:cnf-auto-regressive}, for different times $t_1$ and $t_2$, the outputs of the conditioning networks and their corresponding invertible transformations are different. Therefore the distributions generated at different times are not the same. Moreover, as the coupling is directed, another auto-regressive layer implementing the other coupling direction (from $x_2$ to $x_1$) is needed to guarantee the expressivity. In this work, we apply a network structure that has two auto-regressive layers with reverse coupling orders $(1, 2, ..., n)$ and $(n, n-1, ..., 1)$, so that any two dimensions are coupled. Our model is simpler than the network in \cite{HUANG2023112155} with eight auto-regressive layers, while producing competitive results. Within each auto-regressive layer, we use $5$ bins for each spline transformation layer; the conditioner network contains $2$ hidden layers with width $16$. Consequently, the total number of parameters of the network is $1200$ for all of our $2$D problems. The prior $q$ is chosen to be the standard Gaussian. For all experiments, we use the Adam \cite{DBLP:journals/corr/KingmaB14} optimizer with a learning rate of $10^{-3}$.

% Sampling from distributions $p(\mfx,t)$ at different time $t$ is straight-forward in VCNF as sampling normalizing flow with extra attention to include the temporal dimension in the input.

Our VCNF can also approximate key ingredients of MFC problems efficiently. The velocity field, which is the temporal derivative of the push forward map $\pt f(\mfz,t)$, can be calculated by auto-differentiation or numerical differentiation schemes. Throughout the paper, we use the following central difference schemes to approximate the temporal derivative of the push forward map $f(\mfz, t)$ and the spatial derivative for the log-density $\log p(\mfx, t)$:
\begin{equation}\label{eq:finite_difference}
\begin{aligned}
\pt  f(\mfx, t) & \approx \frac{f(\mfx, t+ \Delta t/2) - f(\mfx, t - \Delta t/2)}{\Delta t} =: \D_{t}^{\Delta t} f(\mfx, t)\,,       \\
\partial_{\mfx_i} \log p(\mfx, t) & \approx \frac{\log p(\mfx+ \Delta x
\mfe_i/2, t) - \log p(\mfx - \Delta x
\mfe_i/2, t)}{\Delta x} =: \lp \D_{\mfx}^{\Delta x} \log p(\mfx, t) \rp_i\,,\forall i 
\end{aligned}
\end{equation}
where $\mfe_i$ is the standard basis in $\mbR^n$, and $\log p$ is evaluated via \cref{equ:log-den}. Applying chain rule to \cref{equ:auto-regressive}, one obtains $\log \left| \text{det} \frac{\partial f_\psi}{\partial \mfz} \right| = \sum_{i=1}^l \log \left| \text{det} \frac{\partial f_\psi^{(i)}}{\partial \mfz^{(i)}} \right|$ where $f_\psi^{(i)}$ is the i-th auto-regressive layer and $\mfz^{(i)} = f^{(i-1)}_\psi \circ \cdots \circ f^{(1)}_\psi(\mfz,t)$ (with convention $\mfz^{(1)}=\mfz$) is the spatial input of $f^{(i)}_\psi$.

\subsection{Numerical approximation for mean field control problems}\label{sec:VCNF quantities}
In this section, we explain the numerical approximation to the MFC problems. We use the dynamical formulation of OT as an illustrative example.

As per \eqref{eq:OT}, we need to calculate the kinetic energy associated with the temporal normalizing flow and enforce the constraints. From the Lagrangian perspective, the continuity equation is inherently satisfied by the flow model. However, the initial and terminal conditions at $t=0,1$ must be enforced explicitly. To address this, we introduce two extra terms of KL divergence to penalize the deviation of the VCNF model at $t=0$ and $1$. The objective functional we used is written in time-continuous form as
\begin{equation*}
\int_0^1  \int \half \norml\bm{v}(\mfx, t)\normr^2 p(\bm{x}, t)\rd\bm{x} \,\rd t + \lambda \lp \KL( p_0 || p(\cdot, 0)) + \KL( p_1 || p(\cdot, 1))\rp\,,
\end{equation*}
where $\lambda>0$ is a weight parameter. While one can use primal-dual algorithms to update the dual variable $\lambda$, we find in experiments that using a fixed $\lambda$ yields better performance. Note that when implementing the KL divergence numerically, the term $\int \log(p_0(\mfx)) p_0(\mfx) \rd x$ in $\KL( p_0 || p(\cdot, 0)) = \int \log(p_0(\mfx)) p_0(\mfx) \rd x - \int \log(p(\mfx,0)) p_0(\mfx) \rd x$ can be omitted as it does not depend on the trainable parameters.

Next, as the kinetic energy contains double integration over both the spatial and temporal dimensions, we employ a Monte Carlo sampling method. The integration over $t$ is simulated through uniform sampling from $[0, 1]$. As for the integration over $\mfx$, we first sample $\{\mfz_j\}_{j=1}^{N_k}$ from latent distribution $q(\cdot)$ and then compute the transformed samples using the normalizing flow model $\mfx=f(\mfz, t)$, as described in \eqref{eq:MFC_Lagrange}. 

Lastly, while the velocity field can be computed via auto-differentiation, we find that the finite difference scheme yields satisfactory results while being computationally more efficient.

Combining all components, we compute the loss function in a mini-batch setting through the Monte Carlo sampling method, ensuring a balance between computational efficiency and approximation accuracy. We first obtain i.i.d. samples $\{ \mfz_j^{(i)} \}_{i=1,\,j=1}^{N_t\quad N_k}$ and $\{\mfz_j^{(b,0)}, \mfz_j^{(b,1)} \}_{j=1}^{N_b}$ from $q(\mfz)$, and then approximate the loss function as
\begin{equation}\label{eq:ot-loss}
\begin{aligned}
&\frac{1}{N_tN_k}\sum_{i=1}^{N_t}\sum_{j=1}^{N_k} \half \norml \D_{t}^{\Delta t} f(\mfz_j^{(i)},
t^{(i)}) \normr_2^2 - \frac{\lambda}{N_b} \sum_{j = 1}^{N_b} \lp\log
p(\mfx_j^{(0)}, 0) + \log p(\mfx_j^{(1)}, 1)\rp \,,   
\end{aligned}   
\end{equation}
where $\mfx_j^{(0)} = f(\mfz_j^{(b, 0)},0)$ and $\mfx_j^{(1)} = f(\mfz_j^{(b, 1)},1)$.
Here, $N_k$ is the batch size of the samples from the normalizing flow models, while $N_b$ is the batch size of the samples from the initial and target distributions. In our experiments, we set $N_t = 20, N_k = 64, N_b= 2048$.

\subsection{Algorithm}
With the above discussions, we present \cref{alg:pseudocode} to solve the MFC problem using VCNF. We illustrate our algorithm using the OT problem. Other related MFC problems can be treated similarly by modifying the loss function in step 6.
\begin{algorithm}
\caption{Solving MFC with VCNF}
\begin{algorithmic}[1]
\STATE{\textbf{Input}: time horizon $T$; time batch size $N_t$; spatial batch size $N_k$, $N_b$; number of iterations $N$; penalty strength $\lambda$; problem dependent parameters such as $\gamma$. }
\FOR{$k=1,2,\ldots, N$}
\STATE{Sample $t^{(1)}, ..., t^{(N_t)}$ uniformly from $[0, T]$.}
\STATE{For each $t^{(i)}$, generate $N_k$ samples $\{\mfz_j^{(i)}\}_{j=1}^{N_k}$ from reference measure $q(\mfz)$.}
\STATE{For $t = 0, 1$, feed the samples $\{\mfz^{(b, 0)}_j, \mfz^{(b, 1)}_j\}_{j=1}^{N_b}$ through the VCNF to obtain samples $\{\mfx_j^{(0)}, \mfx_j^{(1)}\}_{j=1}^{N_b}$ from $p(\cdot|t=0), p(\cdot|t=1)$.}
\STATE{Evaluate the loss function through \cref{eq:ot-loss}}
\STATE{Update the parameters in VCNF using Adam algorithm }
\ENDFOR
\STATE{\textbf{Output}: Optimized VCNF for solving the MFC problem}
\end{algorithmic}
\label{alg:pseudocode}
\end{algorithm}

\section{Numerical experiments}\label{sec:numerics}
In this section, we perform numerical experiments to test the effectiveness of our VCNF. All the experiments were conducted on a single NVIDIA A100 GPU with 40GB memory and CUDA 12.2. with driver 535.129.03. For all experiments, we use the Adam \cite{DBLP:journals/corr/KingmaB14} optimizer with a learning rate of $10^{-3}$ for $30000$ optimization steps. The code reproducing all the experiments can be found in \cite{cnf_ot}.

\subsection{Optimal transport}
We present the numerical results for OT problems in this subsection. We first test the OT problems between two Gaussian distributions, where a closed-form solution is used to compute the errors. We then test more general distributions and present qualitative results.

\medskip
\noindent\textbf{2D OT: Gaussian to Gaussian.} Given two Gaussian distributions $\mcN(\mu_0,\Sigma_0)$ and $\mcN(\mu_1,\Sigma_1)$, the optimal transport map is given by
$$T(x) = \mu_1 + A(\mfx-\mu_0)\,, \quad A = \Sigma_0^{-1 / 2}\left(\Sigma_0^{1 / 2} \Sigma_1 \Sigma_0^{1 / 2}\right)^{1 / 2} \Sigma_0^{-1 / 2}.$$
The corresponding optimal flow map is $f(\cdot,t) = (1-t) \text{Id} + t T(\cdot)$. The squared Wasserstein-$2$ distance between $\mcN(\mu_0,\Sigma_0)$ and $\mcN(\mu_1,\Sigma_1)$ is 
$$\left\|\mu_0-\mu_1\right\|_2^2+\operatorname{Tr}\left(\Sigma_0+\Sigma_1-2\left(\Sigma_0^{1 / 2} \Sigma_1 \Sigma_0^{1 / 2}\right)^{1 / 2}\right).$$
We calculate the accuracy of our method using the above formulations for several different initial and target distributions listed in \cref{table:accuracy-2d-ot}. The benchmarks are squared Wasserstein-2 distances between the initial and target distributions multiplied by $\frac{1}{2}$. The approximated value of squared Wasserstein-2 distances for different penalty $\lambda$ and the relative errors (in bracket) are presented in the table. We leave the errors for the transportation maps in the appendix. 

\begin{figure}[ht]
\centering
\includegraphics[width=\textwidth]{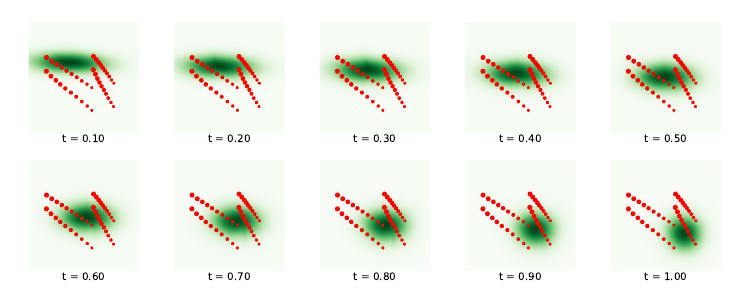}
\caption{2D Gaussian-to-Gaussian in OT problem. The source and target
distributions are given by $\mcN\lp\begin{pmatrix}
-3 \\ 3
\end{pmatrix}, \begin{pmatrix}
5 & 1 \\ 1 & 0.5
\end{pmatrix}\rp$ and $\mcN\lp\begin{pmatrix}
0 \\ 0
\end{pmatrix}, \begin{pmatrix}
1 & 0 \\ 0 & 1
\end{pmatrix}\rp$ respectively. We visualize both the density evolution and several trajectories of the particles. In all the figures, the density heatmap is plotted for the domain $[-6, 3] \times [-3, 6]$. The four red curves represent the trajectories of four starting points $(-5, 3.5), (-5, 2.5), (-1, 3.5), (-1, 2.5)$ moving from upper left to lower right.}
\label{fig:g2g_2D_ot}
\end{figure}

\begin{table}[ht]
\centering
% \small
\normalsize
\resizebox{0.95\columnwidth}{!}{
\begin{tabular}{p{1cm} p{2.5cm} p{2.5cm} p{2.5cm} p{2.5cm} p{2.5cm} p{2.6cm} p{2.6cm}}
\toprule [1.5pt]
\parbox{1cm}{ } &   \parbox{2cm}{ \centering Benchmark}   & \parbox{2.5cm}{ \centering $\lambda=100$ } & \parbox{2.5cm}{ \centering $\lambda=200$} & \parbox{2.5cm}{ \centering $\lambda=500$} & \parbox{2.5cm}{\centering  $\lambda=1000$ } & \parbox{2.6cm}{\centering  $\lambda=2000$ } & \parbox{2.6cm}{\centering  $\lambda=5000$ }  \\ \midrule[1.5pt]
\parbox{1cm}{case 1} 
&  36.000 &  $34.594 (3.90\%)$ &  $35.296 (1.96\%)$ & $35.868 (0.37\%)$ & $35.882 (0.33\%)$ & $36.341 (0.95\%)$ & $36.568 (1.58\%)$   \\ \midrule[0.5pt]
\parbox{1cm}{case 2} 
&  36.125 &  $35.045 (2.99\%)$ &  $35.606 (1.43\%)$ & $36.281 (0.43\%)$ & $36.474 (0.97\%)$ & $36.799 (1.87\%)$ & $37.030 (2.51\%)$    \\ \midrule[0.5pt]
\parbox{1cm}{case 3} 
& 36.860 &  $35.133 (4.69\%)$ &  $35.572 (3.49\%)$ & $36.226 (1.72\%)$ & $35.898 (2.61\%)$ & $36.892 (0.09\%)$ & $37.871 (2.75\%)$   \\ \midrule[0.5pt] 
\parbox{1cm}{case 4} 
& 36.931 &  $33.560 (9.13\%)$ &  $35.370 (4.23\%)$ & $36.555 (1.02\%)$ & $36.700 (0.63\%)$ & $37.090 (0.43\%)$ & $37.380 (1.22\%)$   \\ \midrule[0.5pt]
\parbox{1cm}{case 5} 
& 0.125 &  $0.1254 (0.32\%)$ &  $0.1264 (1.12\%)$ & $0.1284 (2.72\%)$ & $0.1289 (3.12\%)$ & $0.1323 (5.84\%)$ & $0.1347 (7.76\%)$ \\ \midrule[0.5pt]
\parbox{1cm}{case 6} 
& 0.860 &  $0.8345 (2.91\%)$ &  $0.8715 (1.40\%)$ & $0.9068 (5.47\%)$ & $0.9351 (8.72\%)$ & $0.9571 (11.28\%)$ & $0.9973 (15.93\%)$ \\ \midrule[0.5pt]
\parbox{1cm}{case 7} 
& 0.931 &  $0.8907 (4.33\%)$ &  $0.9232 (0.84\%)$ & $0.9564 (2.73\%)$ & $0.9733 (4.54\%)$ & $1.0300 (10.63\%)$ & $1.0779 (15.78\%)$ 
\\\bottomrule[1.5pt]
\end{tabular}
}
\caption{Accuracy of VCNF on 2D Gaussian-to-Gaussian in OT problem. The source and target distributions are Gaussian distributions with different mean and covariances. Mean: in cases 1-4, the source is centered at $(-3,-3)$, and the target is centered at $(3,3)$; in cases 5-7, both the source and the target are centered at $(0,0)$. Covariances: the source covariance for case 1 is $[[1,0],[0,1]]$; the source covariance for case 2 and 5 is $[[1,0],[0,0.25]]$; the source covariance for case 3 and 6 as $[[4,1.5],[1.5,3]]$; the source covariance for case 4 and 7: $[[5,1],[1,0.5]]$. The covariances for the target distributions are identity matrices for all cases. The $N_t$ is fixed to be $20$.}
\label{table:accuracy-2d-ot}
\end{table}
\cref{fig:g2g_2D_ot} shows a visualization of the optimal transportation process. The red dots show the particle trajectories from the source distribution to the target distribution. The densities are also visualized at each time stamp. Greater dots represent earlier time stamps and this convention is kept for later visualization.

% \begin{figure}[h]
%     \centering
%     \includegraphics[width=\textwidth]{fig/velo.pdf}
%     \caption{2D Gaussian-to-Gaussian OT velocity field}
%     \label{fig:g2g_2D_ot_velocity}
% \end{figure}

\begin{figure}[ht]
\centering
\includegraphics[width=\textwidth]{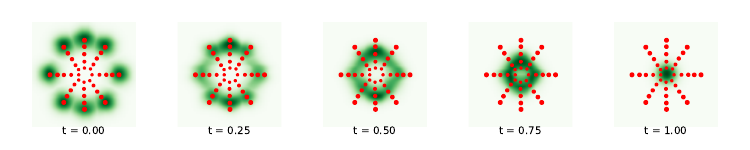}
\caption{2D Gaussian mixture-to-Gaussian in OT problem. All trajectories start from the centers of the Gaussian mixture components, i.e. $(5, 0)$, $(3, 4)$, $(0, 5)$, $(-3, 4)$, $(-5, 0)$, $(-3, -4)$, $(0, -5)$, $(3, -4)$ and the plot domain is $[-7.5, 7.5] \times [-7.5, 7.5].$ The terminal Gaussian distribution is centered at the origin and all Gaussian distributions have covariance matrix $\mfI$.}
\label{fig:gm2g_2D_ot_traj}
\end{figure}

\medskip
\noindent\textbf{2D OT: Gaussian mixture to Gaussian.}
VCNF can also solve complex transport problems from a Gaussian mixture distribution to a Gaussian distribution as shown in \cref{fig:gm2g_2D_ot_traj}. We compute the transport map from a mixture of eight Gaussian distributions with centers $(5, 0), (3, 4), (0, 5), (-3, 4), (-5, 0), (-3, -4), (0, -5), (3, -4)$ to a single Gaussian distribution. All eight trajectories in red dots and the density plot at each time stamp demonstrate the effectiveness of our VCNF model.  

\subsection{Regularized Wasserstein proximal operator}
In this subsection, we compute the RWPO problem \cref{equ:rwpo-loss} using VCNF. We take the initial distribution $p_0 \sim \mathcal{N}(0, 2(T+1) I_d / \beta)$ and the potential function $V(\mfx) = \norm{x}^2/2$. To estimate the velocity field and the score function in \eqref{eq:RWPO}, we apply the finite difference scheme \eqref{eq:finite_difference}.
The reason we avoid using auto-differentiation for the score function is that, the number of forward calculations is proportional to the dimension of the state space, which is not efficient in high dimensions. We also remark that gradient estimators based on control variates can mitigate this issue \cite{Grathwohl2017BackpropagationTT}.

Denote by $\psi$ all the trainable parameters. The objective functional \eqref{eq:RWPO} is approximately
\begin{equation}
\begin{aligned}
L(\psi) = & \ \int_0^T \int_{\RR^d} \half\norm{\p_t f_{\psi}
(\mfz, t) + \frac{1}{\beta} \nabla_{\mfx} \log p(f_\psi(\mfz, t)) }^2 q(\mfz)\,\rd\mfz\, \rd t \\
& + \lambda \KL( p_0 || p(\cdot, 0)) +
\mbE_{\mfx \sim f_{\psi}(\cdot, 1)_\#p}[V(\mfx)]      \\
\approx & \ \frac{T}{2N_t N_k}\sum_{i=1}^{N_t} \sum_{j=1}^{N_k} \norml \D_{t}^{\Delta t} f(\mfz_j^{(i)},
t^{(i)}) + \frac{1}{\beta}\D_{\mfx}^{\Delta \mfx} \log p(f_{\psi}
(\mfz_j^{(i)}, t^{(i)}), t^{(i)}) \normr^2      \\
&  + C + \frac{\lambda}{N_b} \sum_{i = 1}^{N_b}(-\log
p(\mfx_i^{(0)}, 0)) + \frac{1}{N_1}\sum_{j=1}^{N_1} V(f_{\psi}(\mfz_j^{(1)}, 1))\,,
\end{aligned}
\end{equation}
where $N_b, N_1, N_k$ stands for the batch size of the initial condition loss, potential loss, and kinetic energy loss of the objective function, respectively. The constant $C$ represents the entropy of $p_0$ which is irrelavant to the optimization.

\begin{figure}[h]
\centering
\includegraphics[width=\textwidth]{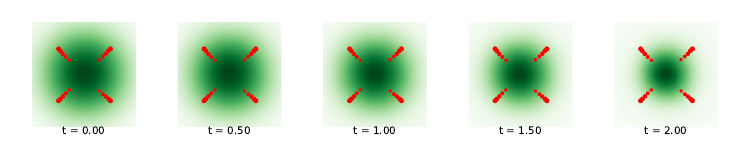}
\caption{RWPO with initial condition of Gaussian distribution with covariance $\Sigma = 4\mfI$ and quadratic potential $V(\mfx) = \| \mfx \|_2^2/2$. We choose $T = 2, \beta = 1, \lambda = 200$ with t batch size $20$. The density on the domain $[-4, 4] \times [-4, 4]$ is visualized with trajectories starting from $(-3, -3), (-3, 3), (3, 3), (3, -3)$ that contract to the middle.}
\label{fig:g2g_2D_rwpo_evolve}
\end{figure}
We perform numerical experiments with different diffusive constants $\beta$ and report the accuracy in \cref{table:accuracy-2d-rwpo}. We observe that the VCNF model gives reasonable results for a wide range of the penalty parameter $\lambda$. 
\begin{table}[tb]
\centering
% \small
\resizebox{0.95\columnwidth}{!}{
\begin{tabular}{p{3cm} p{2cm} p{2.5cm} p{2.5cm} p{2.5cm} p{2.6cm} p{2.5cm}}
\toprule [1.5pt]
\parbox{3cm}{  } &   \parbox{2cm}{ \centering Benchmark}   & \parbox{2cm}{ \centering $\lambda=100$ } & \parbox{2cm}{ \centering $\lambda=200$} & \parbox{2cm}{ \centering $\lambda=500$} & \parbox{2.6cm}{\centering  $\lambda=1000$ } & \parbox{2cm}{\centering  $\lambda=2000$ } % & \parbox{2cm}{\centering  $\lambda=5000$ }
\\ \midrule[1.5pt]
\parbox{3cm}{$\beta = 1, T=1$} 
&  3.386 &  $3.360 (0.79\%)$ &  $3.380 (0.18\%)$ & $3.3984 (0.37\%)$ & $3.404 (0.56\%)$ & $3.460 (2.18\%)$ % & $3.555 (3.99\%)$   
\\ \midrule[0.5pt]
\parbox{3cm}{$\beta = 0.5, T=1$} 
& 6.773 & $7.254 (7.10\%)$ &  $6.856 (1.23\%)$ & $7.570 (11.77\%)$ & $7.016 (3.60\%)$ & $6.872 (1.46\%)$ %& $7.173 (5.92\%)$   
\\ \midrule[0.5pt]
\parbox{3cm}{$\beta = 1, T = 2$}
& 4.197 & $4.197 (0.01\%)$ &  $4.184 (0.31\%)$ & $4.217 (0.46\%)$ & $4.241 (1.05\%)$ & $4.284 (2.06\%)$ %& $4.300 (2.45\%)$   
\\ \midrule[0.5pt]
\parbox{3cm}{$\beta = 0.5, T = 2$}
& 8.394 & $8.944 (6.54\%)$ &  $8.409 (0.18\%)$ & $9.025 (7.51\%)$ & $10.207 (21.59\%)$ & $9.559 (13.88\%)$ %& $18.354 (118.6\%)$  
\\\bottomrule[1.5pt]
\end{tabular}
}
\caption{Accuracy of VCNF on solving the RWPO with a quadratic potential. $T$ is the time period and $\beta$ is the inverse diffusion coefficient. The $N_t$ is fixed to be $20$. It can be observed from the table that as $T$ increases and $\beta$ decreases, the accuracy of our model deteriorates. This may be caused by the density being subject to greater evolutions.}
\label{table:accuracy-2d-rwpo}
\end{table}
% \begin{table}[tb]
% \caption{We evaluate the accuracy of our model in solving the RWPO over different settings. For all the experiments, we use learning rate $10^{-3}$ and epochs $30000$. \XZ{Compare with previous table.} }
% \centering
% % \small
% \resizebox{0.95\columnwidth}{!}{
% \begin{tabular}{p{3cm} p{2.5cm} p{2.5cm} p{2.5cm} p{2.5cm} p{2.5cm} p{2.5cm} p{2.5cm}}
% \toprule [1.5pt]
% \parbox{3cm}{  } &   \parbox{2cm}{ \centering Benchmark}   & \parbox{2cm}{ \centering $\lambda=100$ } & \parbox{2cm}{ \centering $\lambda=200$} & \parbox{2cm}{ \centering $\lambda=500$} & \parbox{2cm}{\centering  $\lambda=1000$ } & \parbox{2cm}{\centering  $\lambda=2000$ } & \parbox{2cm}{\centering  $\lambda=5000$ }  \\ \midrule[1.5pt]
% \parbox{3cm}{$\beta = 1, T=1$} 
%  &  3.386 &  $3.357 (0.87\%)$ &  $3.382 (0.12\%)$  & $3.405 (0.55\%)$ & $3.409 (0.67\%)$ & $3.465 (2.32\%)$ & $3.426 (1.16\%)$   \\ \midrule[0.5pt]
% \parbox{3cm}{$\beta = 0.5, T=1$} 
%  & 6.773 & $6.642 (1.93\%)$ &  $6.850 (1.15\%)$ & $\bold{6.790 (0.26\%)}$ & $\boldsymbol{8.906 (31.5\%)}$ & $\bold{6.902 (1.92\%)}$ & $\bold{16.420 (142.4\%)}$   \\ \midrule[0.5pt]
% \parbox{3cm}{$\beta = 1, T = 2$}
%  & 4.197 & $? (0.01\%)$ &  $? (0.31\%)$ & $? (0.46\%)$ & $4.241 (1.05\%)$ & $? (2.06\%)$ & $? (2.45\%)$   \\ \midrule[0.5pt]
%  \parbox{3cm}{$\beta = 0.5, T = 2$}
%  & 8.394 & $? (6.54\%)$ &  $? (0.18\%)$ & $? (205.1\%)$ & $? (21.59\%)$ & $? (957.7\%)$ & $? (118.6\%)$  
% \\\bottomrule[1.5pt]
% \end{tabular}

% }
% \end{table}

\begin{figure}[h]
\centering
\includegraphics[width=\textwidth]{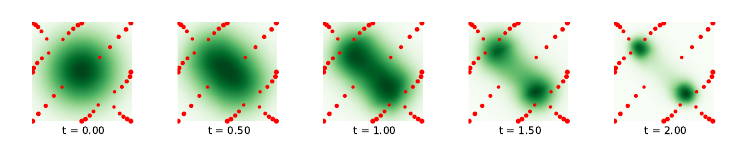}
\caption{RWPO with initial condition of Gaussian distribution with covariance $\Sigma = 4/5\mfI$ and double-well potential $V(\mfx) = ((x_1-1)^2+(x_2+1)^2)((x_1+1)^2 + (x_2-1)^2) / 4$. We choose $T= 2, \beta = 5, \lambda = 100$, and $N_t = 10$. The density on domain $[-2, 2] \times [-2, 2]$ is visualized and red points denote the trajectories starting from the boundary points $(-2, -2), (-2, 0), (-2, 2), (0, 2), (2, 2), (2, 0), (2, -2), (0, -2)$ and then contracting to two wells.}
\label{fig:g2g_2D_rwpo_dw}
\end{figure}

\medskip
\noindent\textbf{Regularized Wasserstein proximal operator with double well potential}. Next, we consider the experiments in \cite{zhou2024score}. We can set $V$ to be a double-well potential, given by $V(\mfx) = ((x_1-a)^2+(x_2+a)^2)((x_1+a)^2 + (x_2-a)^2) / 4$. The explicit solution is stated in \cref{prop:kernel-solutions}. The value of $a$ is taken to be $0.5, 1$. Our solution's visualization and accuracy summary are presented in \cref{fig:g2g_2D_rwpo_dw} and \cref{table:accuracy-2d-rwpo_dw}.

\begin{table}[tb]
\centering
\resizebox{0.95\columnwidth}{!}{
\begin{tabular}{p{3cm} p{2cm} p{2cm} p{2cm} p{2cm} p{2cm} p{2cm}}
\toprule [1.5pt]
\parbox{3cm}{ Density }   & \parbox{2cm}{ \centering $\lambda=100$ } & \parbox{2cm}{ \centering $\lambda=200$} & \parbox{2cm}{ \centering $\lambda=500$} & \parbox{2cm}{\centering  $\lambda=1000$ } & \parbox{2cm}{\centering  $\lambda=2000$ } & \parbox{2cm}{\centering  $\lambda=5000$ }  \\ \midrule[1.5pt]
\parbox{3cm}{$\beta=5, a = 1$} 
& 2.621\% &  2.948\% & 3.067\% & 3.281\% & 3.307\% & 4.545\%   \\ \midrule[0.5pt]  
\parbox{3cm}{$\beta=5, a = 0.5$} 
& 4.129\% &  4.262\% & 4.309\% & 4.357\% & 4.933\% & 4.755\%   \\ \midrule[0.5pt] 
\parbox{3cm}{$\beta=10, a = 1$} 
& 2.766\% &  3.437\% & 4.146\% & 5.566\% & 6.927\% & 5.238\%   \\ \midrule[0.5pt] 
\parbox{3cm}{$\beta=10, a = 0.5$} 
& 4.008\% &  3.858\% & 5.066\% & 4.916\% & 4.472\% & 5.287\%  
\\\bottomrule[1.5pt]
\end{tabular}}
\caption{Accuracy of VCNF on solving RWPO with Gaussian initial condition and double-well potential. Throughout the experiments, we change $a$ (the double-well potential) and $\beta$ (inverse diffusion coefficient) to test the accuracy and robustness of VCNF. $T$ is fixed to be $2$ and $N_t = 10$.}
\label{table:accuracy-2d-rwpo_dw}
\end{table}

\subsection{Flow matching for Fokker--Planck equation}
Lastly, we compute the flow matching problem for the FP equation using VCNF. With the help of the score function, one can rewrite the FP equation \cref{equ:fp-score} as a continuity equation. Therefore, solving the FP equation is equivalent to matching the velocity field of VCNF to the drift field in the FP equation subtracted by the score function. Note that the score function can be efficiently evaluated using VCNF. Recently, solving FP equations under Lagrangian coordinate with neural network is also investigated in \cite{zuo2024numerical}, with numerical analysis on the convergence order of the proposed method.

\medskip
\noindent\textbf{Fokker--Planck equations from Ornstein--Uhlenbeck (OU) processes.} We first consider the FP equation for the OU process $\p_t  p(\mfx, t) = \nabla_{\mfx} \cdot (a\mfx \, p(\mfx, t)) + \gamma \Delta_\mfx p(\mfx, t)\,, t \in [0, T]\,.$
By rewriting it as a continuity equation with velocity given by $-a\mfx - \gamma
\nabla_\mfx \log  p(\mfx, t)$, we can solve it by flow matching with the loss function
\begin{equation*}
\int_0^T  \int \norml \bm{v}(\mfx, t) + a\mfx
+ \gamma\nabla_{\mfx} \log  p(\mfx, t) \normr^2  p(\bm{x}, t)\,\rd\bm{x}\, \rd t + \lambda \KL( p_0 ||  p(\cdot, 0))\,.
\end{equation*}
We remark that this loss function is different from the physics-informed loss function \cite{zeng2023adaptive}, which is derived from the FP equation via auto-differentiation. Instead, our method only requires first-order derivatives, which is computationally more efficient. This loss is also similar to the flow matching loss function in~\cite{lipman2022flow}. We approximate the integration via Monte Carlo sampling method: 
\begin{equation}\label{equ:fp_loss}
\begin{aligned}
L(\psi)
\approx & \ \frac{1}{N_t N_k}\sum_{i=1}^{N_t} \sum_{j=1}^{N_k} \norml \p_t f_{\psi}(\mfz_j^{(i)}, t^{(i)})
+ a f_{\psi}(\mfz_j^{(i)}, t^{(i)}) + \gamma\nabla_{\mfx} \log p(f_{\psi}(\mfz_j^{(i)}, t^{(i)})
, t^{(i)})  \normr^2      \\
& \ + \lambda \sum_{i = 1}^{N_b} (-\log
p(\mfx_i^{(0)}, 0))\,.
\end{aligned}
\end{equation}
The score function is estimated via finite differences \eqref{eq:finite_difference}. We report the root mean squared error (RMSE) of the solution of FP equation using our model at time $T = 1$ over a grid of size $N = 500 \times 500$ on $[-5, 5] \times [-5, 5]$. The RMSE error between the computed density $p(\cdot,t) = f_\psi(\cdot,t)_\# q(\cdot)$ and the true solution $p^*$ at $t=1$ is approximated by Monte Carlo sampling:
\begin{equation*}
\text{RMSE}(p^*, p) \approx \sqrt{\frac{1}{N}\sum_{i=1}^N\lp p^*(\mfx_i, t=1) -  p(\mfx_i, t=1)\rp^2}\,,
\end{equation*}
where $\mfx_i = f_\psi(\mfz_i,1)$ and $\mfz_i$ are i.i.d. sampled from $q(\cdot)$.
The results are concluded in \cref{table:accuracy-2d-fk}. We observe that our model is robust and is not sensitive to the penalty parameter $\lambda$. 
\begin{table}[tb]
\centering
\resizebox{0.95\columnwidth}{!}{
\begin{tabular}{p{2cm} p{2cm} p{2cm} p{2cm} p{2cm} p{2cm} p{2cm}}
\toprule [1.5pt]
\parbox{2cm}{ Density }   & \parbox{2cm}{ \centering $\lambda=100$ } & \parbox{2cm}{ \centering $\lambda=200$} & \parbox{2cm}{ \centering $\lambda=500$} & \parbox{2cm}{\centering  $\lambda=1000$ } & \parbox{2cm}{\centering  $\lambda=2000$ } & \parbox{2cm}{\centering  $\lambda=5000$ }  \\ \midrule[1.5pt]
\parbox{2cm}{$a = 1$} 
& \num{8.274E-4} &  \num{5.966E-4} & \num{5.840E-4} & \num{6.776E-4} & \num{6.521E-4} & \num{1.074E-3}   \\ \midrule[0.5pt]  
\parbox{2cm}{$a = 0.5$} 
& \num{3.750E-4} &  \num{1.482E-4} & \num{2.477E-4} & \num{ 4.138E-4} & \num{4.354E-4} & \num{2.812E-3} 
\\\bottomrule[1.5pt]
\end{tabular}}
\caption{Accuracy of VCNF on solving the FP equation via the flow matching framework. We vary the drift parameter $a$. The initial distribution is Gaussian with zero mean and covariance $4\mfI$ and $\gamma=0.5$. This table records the RMSE of the solution obtained by VCNF calculated on a $500 \times 500$ grid of $[-5, 5] \times [-5, 5].$}
\label{table:accuracy-2d-fk}
\end{table}

\begin{figure}[h]
\centering
\includegraphics[width=\textwidth]{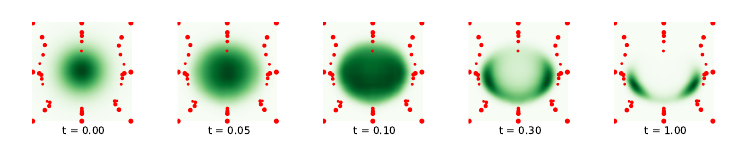}
\caption{FP equation with standard Gaussian initial distribution and ``smiling'' invariant distribution. The ``smiling'' distribution $\pi(\mfx)$ has potential function $V(\mfx) = \frac{1}{4}(x_1^2 + x_2^2 - 4)^2 + (x_2+1)^2$. The $\delta$ is fixed to be $0.5$. Notice that we choose non-uniform time step for better visualization of the trajectories. The domain of visualization is $[-3, 3]\times[-3, 3]$ with all trajectories starting at the boundary points $(-3, -3),(-3, 0),(-3, 3),(0, 3),(3, 3), (3, 0), (3, -3)$. All the trajectories converge to the lower hemicycle $x_1^2 + x_2^2 = 4.$}
\label{fig:nongrad_fp}
\end{figure}

\medskip
\noindent\textbf{Fokker--Planck equation with non-gradient velocity field.}
In the previous example, the velocity field $a\mfx$ can be viewed as the gradient field of some potential function, which gives the FP equation a gradient flow interpretation in Wasserstein-$2$ metric space. Moreover, the invariant measure satisfies the detailed balance condition at equilibrium. However, FP equations with non-gradient velocity fields are also of importance \cite{onsager1931reciprocal, gao2021random}.
Many FP equations admit non-equilibrium stationary states without detailed balance, this is also closely related
to the recent success in generative modeling using diffusion models~\cite{sohl2015deep}. In this section, we consider FP equations with a non-gradient velocity field and a ``smiling'' invariant measure: 
\begin{equation*}
\pi(\mfx) \propto \exp(-U(\mfx)) = \exp\lb -\frac{1}{4}(x_1^2 + x_2^2 - 4)^2 - (x_2+1)^2 \rb.
\end{equation*}
The non-gradient velocity field is given by a combination of the original
gradient vector field and a small Hamiltonian vector field perturbation of size $\delta\in\mathbb{R}$
\begin{equation*}
    \begin{aligned}
        \mfv(\mfx) = -\nabla U(\mfx) - \delta \begin{pmatrix}
0 & 1 \\
-1 & 0
\end{pmatrix}\nabla U(\mfx) = -\begin{pmatrix}
    (x_1+\delta x_2)(x_1^2 + x_2^2 - 4) + 2\delta (x_2+1) \\
    (x_2 - \delta x_1)(x_1^2 + x_2^2 - 4) + 2(x_2+1)
\end{pmatrix}\,.
    \end{aligned}
\end{equation*}
This Hamiltonian vector field perturbation will leave the original invariant measure unchanged while breaking the detailed balance condition for the stochastic process \cite{gao2021random}. Our results are demonstrated in \cref{fig:nongrad_fp}. Similar to \cref{equ:fp_loss}, the objective function $L(\psi)$ is approximately given by
\begin{equation}
\begin{aligned}
%L(\psi)\approx
& \ \frac{1}{N_t N_k}\sum_{i=1}^{N_t} \sum_{j=1}^{N_k} \norml \p_t f_{\psi}(\mfz_j^{(i)}, t^{(i)})
+ \gamma\nabla_{\mfx} \log p(f_{\psi}(\mfz_j^{(i)}, t^{(i)})
, t^{(i)}) - \mfv(f_{\psi}(\mfz_j^{(i)}, t^{(i)})) \normr^2      \\
& \ + \lambda \sum_{i = 1}^{N_b} (-\log
p(\mfx_i^{(0)}, 0))\,,
\end{aligned}
\end{equation}
with $a(f_{\psi}(\mfz_j^{(i)}, t^{(i)}))$ substituted by $-\mfv(f_{\psi}(\mfz_j^{(i)}, t^{(i)}))$.

% \medskip
% \noindent\textbf{Lorenz system}
% Consider the Lorenz-63 system: 
% \begin{equation*}
% \left\{\begin{aligned}
% \dot x=& c_1(y-x),\\
% \dot y=& x(c_2-z)-y,\\
% \dot z=& xy-c_3 z,
% \end{aligned}\right.
% \end{equation*}
% where we study the parameters as $(c_1, c_2, c_3) = (10, 28, \frac{8}{3})$. In this example, we apply the proposed method to compute the invariant distribution of the above ODE system. \cite{yang2023optimal} has reference solution.

\subsection{Scalability of the algorithm}\label{sec:scalability}
We also test the scalability of the VCNF on solving high dimensional problems. For OT, we can solve for transport maps between Gaussian distributions over 20D within 17 mins to achieve 1\% relative error using 30000 steps. For the FP equation, our model can handle 10D problems in 34 mins for 30000 steps and achieve a $10^{-5}$ absolute $L_2$ error calculated by Monte Carlo sampling with $10^8$ samples. For the RWPO, our algorithm takes 53 mins to run 30000 steps and achieves a relative error of $0.44\%$. We also find that choosing the smallest temporal batch size, i.e., $N_t = 1$ throughout all the high-dimensional experiments is extremely advantageous as it significantly reduces the computational cost while maintaining a good accuracy.

\section{Discussion}\label{sec:discussion} 
In this paper we developed VCNF, which is a neural network-based framework for solving MFC problems using methods in generative models. In particular, our neural network structure makes use of both CNFs and neural spline flows. By leveraging the conditional generative model, one can capture the probability distribution's temporal evolution at any given time, which naturally fits into the trajectory-wise formulation of the MFC problems. Moreover, the CNF structure allows us to evaluate various quantities of interests, such as the velocity, kinetic energy, score function, through numerical differentiation and Monte--Carlo sampling. By taking a Lagrangian perspective of MFC problems, the objective functional can be efficiently estimated and optimized by Monte--Carlo sampling. We demonstrate the effectiveness and accuracy of the VCNF in solving various problems, including OT, RWPOs, and controlling non-gradient vector field FP equations. Specifically, we solve the initial value FP equation under the framework of flow matching technique from generative models. Our model is robust to handle evolution problems between distributions with different modalities, e.g., from a Gaussian mixture distribution to a Gaussian distribution. Moreover, our method can be efficiently scaled to high dimensions with reasonable computational cost. 

In future work, we shall explore using VCNF on flow matching problems to simulate general FP equations, including those with nonlinear drift vector fields. Another interesting direction is to apply the proposed framework to real image data sets for solving time-reversible diffusion models.

\section*{Acknowledgement}
J. Zhao, M. Zhou, X. Zuo, and W. Li are supported by AFOSR YIP award No. FA9550-23-1-0087. W. Li is also supported by NSF DMS-2245097, and NSF RTG: 2038080.

% \newpage
% \input{method with time varying parameters}

% \newpage
% \input{numerics}

% \newpage
% \input{2024_08_31}

% \newpage
% \input{numerical_example_Mo}

\newpage
\appendix

\section{Implementation details}

\subsection{Details of the network architecture}\label{sec:neural_spline_flow}
We briefly introduce the conditional normalizing flow
beyond \cref{sec:cnf} with consistent notations. We refer interested readers to \cite{10.5555/3454287.3454962} for more details. We first discuss how the conditioning network is
used to parametrize the monotonic rational-quadratic spline:
\begin{itemize}
    \item A conditioning neural network $\text{NN}^{(k)}$ takes $x_{1:k}$ and $t$ as input and outputs a vector $\theta_k$ of length $3K - 1$.
    \item Vector $\theta_k$ is partitioned as $\theta_k = [\theta_k^w, \theta_k^h, \theta_k^d]$, where $\theta_k^w, \theta_k^h$ have length $K$, and $\theta_k^d$ has length $K - 1$.
    \item Vectors $\theta_k^w, \theta_k^h$ are passed through an elementwise softmax function separately and then multiplied by $2B$. Two outputs $l^{(x)} = (l_1^{(x)}, ..., l_K^{(x)}), (l_1^{(y)}, ..., l_K^{(y)})$ are both vectors with $K$ positive components that sum up to $2B$. Therefore, they are used as the width of $K$ bins spanning $[-B, B]$ in $x$ and $y$ axes, see Figure 1 of \cite{10.5555/3454287.3454962}. Equivalently, the nodes of the spline function are given by
    \begin{equation}
        x_k = \sum_{i=1}^k l_i^{(x)} - B, \quad y_k = \sum_{i=1}^k l_i^{(y)} - B, \quad x_0 = y_0 = -B, x_K = y_K = B. 
    \end{equation}
  %   \mz{don't understand this part.}
    \item The vector $\theta_k^d$ is passed through a softplus function and is interpreted as the values of the derivatives $\{\delta^{(k)}\}_{k=1}^{K-1}$ at the internal knots.
\end{itemize}
The above method constructs a monotonic, continuously-differentiable,
rational-quadratic spline which passes through the knots, with the given
derivatives at the knots. Defining
$s^{(k)} = \frac{y^{(k+1)} - y^{(k)}}{x^{(k+1)} - x^{(k)}}, \xi(x) = \frac{x - x^{(k)}}{x^{(k+1)} - x^{(k)}} \in [0,1]$, the expression for the rational-quadratic spline over the $k$-th bin is given by
\begin{equation}
    \frac{\alpha^{(k)}(\xi)}{\beta^{(k)}(\xi)} = y^{(k)} + \frac{(y^{(k+1)} - y^{(k)})[s^{(k)}\xi^2 + \delta^{(k)}\xi(1-\xi)]}
    {s^{(k)} + [\delta^{(k + 1)} + \delta^{(k)} - 2s^{(k)}]\xi(1-\xi)}
\end{equation}
for $\xi \in [0,1]$. In practice, we use $K=5$ bins for each spline transformation layer.

\section{Mathematical formulation of the numerical experiments}
In this section, we provide the details of the mathematical formulations of the numerical experiments discussed in \cref{sec:numerics}.
\subsection{Regularized Wasserstein proximal operators}\label{sec:rwpo}
Let the initial distribution be $ p_0 \sim \mathcal{N}(0, 2(T+1) I_d / \beta)$. The FP-HJB system governing the RWPO problem is
\begin{equation*}
\left\{\begin{aligned}
&\pt p(\mfx,t) + \nabla_{\mfx}\cdot \parentheses{ p(\mfx,t) \nabla_{\mfx}
\phi(\mfx,t)}=\frac1\beta \Delta_{\mfx}  p(\mfx,t)\,,\\
&\pt\phi(\mfx,t) + \frac12 \norm{\nabla_{\mfx} \phi(\mfx,t)}^2 + \frac1\beta
\Delta_{\mfx}\phi(\mfx,t) =0\,,\\
&  p(0,\mfx)= p_0(\mfx)\,,\quad \phi(\mfx,t)=-V(\mfx)=-\norml \mfx \normr^2/2\,. 
\end{aligned}\right.
\end{equation*}
The true solution $p^*$ is given by
\begin{equation*}
p^*(\mfx,t) = \parentheses{4\pi (T-t+1)/\beta}^{-\frac{d}{2}}
\exp\parentheses{-\dfrac{\beta\norml \mfx \normr^2}{4 (T-t+1)}}\,.
\end{equation*}
The solution to the HJB equation is 
\begin{equation*}
\phi(\mfx,t) = \frac1\beta d\log\frac{1}{T-t+1}-\frac{\norml \mfx \normr^2}{2(T-t+1)}\,.
\end{equation*}
The optimal velocity is $v^* = \nabla \phi$.
Given the close form solution, we can obtain the exact value of the objective functional as
\begin{equation}
\begin{aligned}
& \ \int_0^T \int_{\RR^d} \frac12\norm{v^*(\mfx,t)}^2  p^*(\mfx,t)
\,\rd \mfx \,\rd t + \int_{\RR^d} V(\mfx)  p^*(\mfx,t) \,\rd \mfx     \\
= & \ \int_0^T \int_{\RR^d} \frac{\norml \mfx \normr^2}{2(T-t+1)^2}
 p^*(\mfx,t) \,\rd \mfx \,\rd t + \int_{\RR^d} \frac12
\norml \mfx \normr^2  p^*(\mfx,t) \,\rd \mfx   \\
= & \ \frac{d}{\beta}\lp\int_0^T \frac{\,\rd t}{(T-t+1)} + 1\rp =
\frac{d}{\beta}\lp\log (T + 1) + 1\rp\,.
\end{aligned}
\end{equation}

To solve the RWPO problem using our framework, note that the continuity equation contains a diffusion term, which is equivalent to substracting a score function from the original velocity field 
\begin{equation}
\pt p(\mfx,t) + \nabla_{\mfx}\cdot \parentheses{ p(\mfx,t)
(\mfv(\mfx,t) - \frac1\beta\nabla \log  p(\mfx,t))} = 0\,.
\end{equation}

Consequently, the loss function can be directly modified to include the score function 
\begin{equation} 
\inf_\mfv \int_0^T \int_{\RR^d} \frac12\norm{\p_t f(\mfz, t)
+ \frac1\beta\nabla \log  p(\mfx,t)}^2  p(\mfx,t) \,\rd \mfx \,\rd t + \int_{\RR^d} V(\mfx)  p(\mfx,T) \,\rd \mfx\,.
\end{equation}

For generalized potential function $V(\mfx)$, one has the following kernel solution \cite{li2023kernel}.
\begin{proposition}[Kernel solutions]\label{prop:kernel-solutions}
The solution to the FP-HJB system for RWPO problem is given by the following kernel formulation.
\begin{equation*}
p(\mfx,t)=\frac{1}{(4\pi\frac{1}{\beta}\frac{t(T-t)}{T})^{\frac{d}{2}}}
\int_{\mathbb{R}^d}\int_{\mathbb{R}^d}\frac{e^{-\frac{\beta}{2}(V(\mfy')+
\frac{\|\mfx-\mfy'\|^2}{2(T-t)}+\frac{\|\mfx-\mfy\|^2}{2t})}}{\int_{\mathbb{R}^d}
e^{-\frac{\beta}{2}(V(\tilde \mfy)+\frac{\|\mfy-\tilde \mfy\|^2}{2T})}\rd
\tilde \mfy} p(0,\mfy) \rd \mfy\rd \mfy'\,,
\end{equation*}
and
\begin{equation*}
\phi(\mfx,t)=\frac{2}{\beta} \log\Big(\int_{\mathbb{R}^d}\frac{1}
{(\frac{4\pi}{\beta} (T-t))^{\frac{d}{2}}}e^{-\frac{\beta}{2}(V(\mfy)+
\frac{\|\mfx - \mfy\|^2}{2(T-t)})}\rd \mfy\Big)\,.
\end{equation*}
In particular, 
\begin{equation*}
p(\mfx,T)=\int_{\mathbb{R}^d}K(\mfx, \mfy, \beta, T, V) p(0, \mfy)\rd \mfy\,,
\end{equation*}
where $K\colon \mathbb{R}^d\times\mathbb{R}^d\times \mathbb{R}_+\times \mathbb{R}_+\times C^1(\mathbb{R}^d)\rightarrow\mathbb{R}$ is the kernel function 
\begin{equation*}
K(\mfx,\mfy,\beta, T,V):=\frac{e^{-\frac{\beta}{2}(V(\mfx)+
\frac{\|\mfx-\mfy\|^2}{2T})}}{\int_{\mathbb{R}^d}e^{
-\frac{\beta}{2}(V(\wtd \mfy)+\frac{\|\wtd \mfy - \mfy\|^2}{2T})}\rd \wtd \mfy}\,.
\end{equation*}
The optimal cost is given by $-\int_{\RR^d} \phi(\mfx,0) p(\mfx,0) \,\rd \mfx$.
\end{proposition}

\subsection{Fokker--Planck equation}
We first consider the FP equation corresponding to the Ornstein–Uhlenbeck (OU) process with constant drift $a> 0$ and constant diffusion $\gamma >0$ 
\begin{equation}
\begin{aligned}
\rd\mfx_t & = -a\mfx_t \rd t + \sqrt{2\gamma} \,\rd\mfB_t\,,      \\
\p_t  p(\mfx, t) & = \nx \cdot (a\mfx  p(\mfx, t)) +
\gamma \Delta_\mfx p(\mfx, t)\,.
\end{aligned}
\end{equation}
%\mz{Again, we shall put these details in the appendix}
Using It\^o's calculus, we have
\begin{equation}
\rd\norml \mfx_t \normr^2 = 2(\gamma - a\norml \mfx_t \normr^2) \,\rd t + 2\sqrt{2\gamma} \mfx_t \,\rd\mfB_t.
\end{equation}
Taking expectations on both sides, we obtain that the second moment of $\mfx_t$ satisfies
the following ODE
\begin{equation}\label{eq:OU_ode}
\frac{\rd\mbE \norml \mfx_t \normr^2}{\rd t} = 2(\gamma - a\mbE \norml \mfx_t \normr^2)\,.
\end{equation}
The analytic solution of \cref{eq:OU_ode} is given by $\mbE
\norml \mfx_t \normr^2 = \frac{\gamma}{a} + (\mbE \norml \mfx_0 \normr^2 - \frac{\gamma}{a})e^{-2at}$, which is used to calculate the accuracy of VCNF. 

% \subsection{Potential ablation study on network architecture}
% In this section, we suggests several ablation studies of the network model on
% the numerical experiments:
% \begin{itemize}
%     \item Number of autoregressive layer
%     \item Number of bins: number of bins used to construct the rational quadratic spline;
%     \item Number of MLP conditioner layer: 
%     \item Hidden dimension of MLP conditioner 
% \end{itemize}

\bibliographystyle{siamplain}
\bibliography{references}
\end{document}